# 重构谐波平衡法及其求解复杂非线性问题应用 [1]

代洪华 [2]　王其偲　严子朴　岳晓奎

（西北工业大学航天学院，西安 710072)

（航天飞行动力学技术国家重点实验室，西安 710072）

**摘要**　谐波平衡法是求解非线性动力学系统周期解的最常用方法，但对非线性项进行高阶近似需要庞杂的公式推导，限制了该方法的超高精度解算。作者团队通过对频域非线性量的时域等价重构，提出了重构谐波平衡法（RHB 法），解决了经典谐波平衡法超高阶次计算难题。然而，两种方法均要求动力学系统为多项式型非线性，且无法直接用来求解非线性系统的拟周期解。针对上述问题，本文提出一种将 RHB 法和复杂非线性系统等价重铸法相结合的计算方法，首先将一般非线性问题无损重铸为多项式型非线性系统，然后用 RHB 法进行高精度求解；针对拟周期响应求解问题，提出基于"补频"思想的 RHB 方法，通过基频的优化筛选，实现拟周期响应的快速精准捕捉。选取非线性单摆、相对论谐振子、非线性耦合非对称摆等典型系统进行仿真计算，仿真结果表明，所提出的 RHB-重铸法解非多项式型非线性系统的精度为 $10^{-12}$ 量级，达计算机精度，远超现有方法水平。补频 RHB 法则实现了拟周期问题的高效解算。

**关键词**　非多项式型非线性系统，重构谐波平衡法，微分方程重铸，非线性单摆，拟周期响应

# RECONSTRUCTION HARMONIC BALANCE METHOD AND ITS APPLICATION IN SOLVING COMPLEX NONLINEAR DYNAMICAL SYSTEMS[1]

Dai Honghua[2]　　Wang Qisi　　Yan Zipu　　Yue Xiaokui

( *School of Astronautics*，*Northwestern Polytechnical University*，*Xi'an* 710072，*China* )

( *The National Key Laboratory of Aerospace Flight Dynamics*，*Xi'an* 710072，*China* )

**Abstract**　The harmonic balance method is the most commonly used method for solving periodic solutions of nonlinear dynamic systems, but the high-order approximation of nonlinear terms requires sophisticated formula derivation, which limits its ultra-high accuracy. The authors' team proposed the reconstruction harmonic balance （RHB） method through the equivalent reconstruction of the frequency domain nonlinear quantity in the time domain, which settled the problem of ultra-high-order calculation of the classical harmonic balance method. However, both methods require the dynamical system to be polynomial nonlinear, and cannot be directly used to solve the quasi-periodic solution of the nonlinear system. In view of the above problems, this paper proposes a computational method that combines the RHB method and the recast technique for complex nonlinear systems.

First, the general nonlinear problem is non-destructively recast into a polynomial nonlinear system, and then the RHB method is used for high-precision solutions. Aiming at computing the quasi-periodic response, the RHB method based on the idea of "supplemental frequency" is derived. By optimizing and selecting base frequencies, the fast and accurate capture of quasi-periodic response is achieved. The typical systems such as nonlinear pendulum, relativistic harmonic oscillator, and nonlinear coupling asymmetric pendulum are selected for simulation. The simulation results show that the accuracy of the RHB-recast method for solving nonpolynomial nonlinear systems is on the order of $10^{-12}$, reaching the computer accuracy, far exceeding state-of-the-art methods. The supplemental frequency RHB method achieves the efficient solution of quasi-periodic problems.

**Key words**　nonpolynomial nonlinear system, reconstruction harmonic balance method, differential equation recast technique, nonlinear pendulum, quasi-periodic response

# 引　言

工程中多数动力学问题，其数学模型都是非线性的，线性系统只是在一定假设及限制条件下对非线性系统的理想化近似[1]。随着控制科学与航空宇航任务日益复杂，强非线性的影响越发不容忽视。谐波平衡法（harmonic balance method, HB）在求解非线性系统周期解中应用广泛，但作为一种半解析半数值方法，随着系统自由度、方法阶数的增加，公式推导工作将变得困难[2]。Hall 等[3]提出了高维谐波平衡法（high-dimensional harmonic balance, HDHB），通过"频域非线性量的时域近似计算"来简化公式推导，但是，由于引入近似关系导致非物理假解问题[4-5]。Dai 等[6]发现了频域和时域非线性项之间的条件等价恒等式，基于此，首创了重构谐波平衡法（reconstruction harmonic balance, RHB）实现超高阶（$N>100$）高精计算，并给出时域配点计算的最优采样定理，从理论上消除非物理解。

HB 类方法（HB 法及其改进方法）不仅用于计算简单非线性系统（杜芬、范德波尔方程等），还发展到航空航天、深空探测等前沿领域。哈尔滨工业大学陈毅等[7]使用谐波平衡-交变频域/时域法（HB-alternating frequency-time, HB-AFT）用于求解航空发动机中双转子-轴承-机匣系统动力学方程。中山大学陈衍茂[8]等使用增量谐波平衡法对带机外挂载的二元机翼进行动力学特性研究。RHB 法也被作为三体轨道的设计依据，计算结果符合鹊桥中继卫星实际飞行数据[6]。

但是该方法仍面临着两个问题：1）谐波平衡法的本质依赖于非线性项的傅里叶级数展开，因此，受限于多项式型非线性系统求解。对于非多项式型复杂非线性问题，谐波平衡法难以适用。2）已有谐波平衡类方法建立在单基频的假设上进行级数展开，由于拟周期响应存在多基频的特征，因此已有方法难以直接求解高精度拟周期解。

针对非多项式型非线性系统的求解问题，目前有两类处理方法，分别为直接法和间接法。直接法包括 HB-Taylor 法[9]与 HB-AFT 法[10-11]。HB-Taylor 法通过泰勒级数展开将非线性函数用有限阶近似多项式描述；HB-AFT 法通过对非线性项的时域值采样以离散原问题。由于直接法对原系统进行了近似，导致求解精度低，且计算性能分别受制于高阶级数描述和过采样等问题。Cochelin 等[12]提出的谐波平衡-重铸法（HB-recast）是一种间接方法，通过重铸（recast）技术，成功将复杂非线性微分动力学系统无损变换为多项式型微分代数方程，然后用 HB 法加以求解[13-14]。但是，受限于 HB 法的高阶计算（使用重铸法会增加系统的维度，进一步增加了高阶计算的难度）与原重铸形式的二次型限制（原方法要求新系统中非线性至多为二次多项式），至今难以对复杂非线性系统周期响应进行高效高精求解。

第二个难题是拟周期响应求解问题。拟周期响应大量出现在非线性动力学系统中[15-17]，其频率响应由多个不可约基频及其线性组合描述[18]。由于传统 HB 类方法基于单个基频进行近似解逼近，不能简单通过基频及其整数倍频率分量对拟周期响应描述，Chua 等[19]提出了结合广义傅里叶级数的改进多频 HB 法，实现拟周期响应的准确求解。然而，使用多频 HB 类方法对强非线性项进行频域内高阶描述时，计算效率严重受限于高维傅里叶分析（频域分量由多重积分[20]、求和[21]计算得到）。Liu 等[22]使用多频 HDHB 法求解受迫范德波尔振子的稳态响应，避免非线性项谐波系数的直接表示以提高求

解效率，但是多频 HDHB 法的精度受到非物理解的破坏[23]。总之，当前关于 HB 类方法的研究已拓宽到对拟周期响应求解领域，然而由于多基频计算中的高阶频域描述困难，对此类复杂响应的高性能求解尚待解决。

针对上述问题，本文提出了 RHB 法与重铸法、多频谐波平衡计算相结合的两种方法：（1）RHB-重铸法和（2）重构多谐波平衡法（reconstruction multiple harmonic balance, RMHB）。一方面 RHB-重铸法通过将一般非线性等价转化为多项式型非线性系统，再采用 RHB 法以确定时域最优采样点数，实现对复杂非线性动力学系统的高阶高精求解，仿真误差达计算机精度。另一方面，RMHB 法通过筛选和补充多个基频，借鉴 RHB 法的时域等价重构思想，完善了学术界在使用多频 HB 类方法求解拟周期响应时消除混淆假解的理论研究。

# 1 谐波平衡法及重构谐波平衡法

## 1.1 谐波平衡法

对非线性动力学系统
$$\dot{x} = f(x,t) \tag{1}$$
HB 法用有限阶傅里叶级数来构造近似解及其导数
$$x(t) = x_0 + \sum_{n=1}^{N}[x_{2n-1}\cos n\omega t + x_{2n}\sin n\omega t] \tag{2}$$
$$\dot{x}(t) = \sum_{n=1}^{N}[-n\omega x_{2n-1}\sin n\omega t + n\omega x_{2n}\cos n\omega t] \tag{3}$$
其中 $N$ 是 HB 法的阶数。$x_0, x_1, \cdots, x_{2n}$ 为未知傅里叶系数，也被称为频域变量，将频域变量组成待求解向量 $\hat{x} = [x_0\ x_1\ \cdots\ x_{2n}]^T$。假设多项式函数 $f(x) = x^\phi$，忽略由非线性项而出现的高次谐波，只需展开前 $N$ 次的傅里叶分量：
$$f(x,t) = x^\phi(t) = r_0 + \sum_{n=1}^{N}[r_{2n-1}\cos n\omega t + r_{2n}\sin n\omega t] \tag{4}$$
其中各次分量 $r_0, r_1, \cdots, r_{2n}$ 为

$$\begin{cases} r_0 = \frac{1}{2\pi}\int_0^{2\pi}\left[x_0 + \sum_{n=1}^{N}(x_{2n-1}\cos n\theta + x_{2n}\sin n\theta)\right]^\phi d\theta \\ r_{2n-1} = \frac{1}{\pi}\int_0^{2\pi}\left[x_0 + \sum_{n=1}^{N}(x_{2n-1}\cos n\theta + x_{2n}\sin n\theta)\right]^\phi \cos(n\theta)d\theta \\ r_{2n} = \frac{1}{\pi}\int_0^{2\pi}\left[x_0 + \sum_{n=1}^{N}(x_{2n-1}\cos n\theta + x_{2n}\sin n\theta)\right]^\phi \sin(n\theta)d\theta \end{cases} \tag{5}$$

其中 $n=1, 2, \ldots, N$；$\theta = \omega t$。将表示非线性函数的傅里叶分量记为 $\hat{f}$，是待求解向量 $\hat{x}$ 的非线性函数。

将式(2)-(4)代入微分方程(1)，令常数项及前 $n$ 次谐波 $\cos n\omega t$、$\sin n\omega t (n=1,2,\cdots,N)$ 系数配平，得到非线性代数方程组
$$\omega A\hat{x} = \hat{f}(\hat{x}) \tag{6}$$
其中，分块矩阵 $A$ 为
$$A = \text{dig}([0, J_1, J_2, \cdots, J_N]),\quad J_n = n\begin{bmatrix} 0 & 1 \\ -1 & 0 \end{bmatrix}$$
由三角函数的求导关系得到。

非线性项的近似表示(4)中，符号运算的复杂度会随算法阶次的提高呈指数增长[6]。当 HB 法的阶数超过 20，即便有计算机数学软件的辅助，非线性项的推导整理工作量仍难以接受。

## 1.2 高维谐波平衡法

为简化 HB 法的符号运算量，Hall 等[3]使用时域值替代频域分量，即将 $N$ 阶 HB 法中的傅里叶系数作用转换矩阵 $E_{\text{HDHB}}$，建立与一个周期 $2N+1$ 个等距配点上时域量间的联系，定义
$$\hat{x} = E_{\text{HDHB}}\tilde{x},\quad \hat{f} \approx E_{\text{HDHB}}\tilde{f} \tag{7}$$
其中
$$\tilde{x} = [x(t_0)\ x(t_1)\ \cdots\ x(t_{2N})]^T$$
$$\tilde{f} = [f(x(t_0))\ f(x(t_1))\ \cdots\ f(x(t_{2N}))]^T$$
$$E_{\text{HDHB}} = \frac{1}{2N+1}\begin{bmatrix} 1/2 & 1/2 & \cdots & 1/2 \\ \cos\omega t_0 & \cos\omega t_1 & \cdots & \cos\omega t_{2N} \\ \sin\omega t_0 & \sin\omega t_1 & \cdots & \sin\omega t_{2N} \\ \vdots & \vdots & & \vdots \\ \cos N\omega t_0 & \cos N\omega t_1 & \cdots & \cos N\omega t_{2N} \\ \sin N\omega t_0 & \sin N\omega t_1 & \cdots & \sin N\omega t_{2N} \end{bmatrix}$$

将式(7)代入 HB 法代数方程组(6)得到 HDHB 法求解微分方程(1)对应的代数方程组
$$\omega A E_{\text{HDHB}}\tilde{x} = E_{\text{HDHB}}\tilde{f}(\tilde{x}) \tag{8}$$

HDHB 法显著提高计算效率，并被认为是 HB 法的一种改进，在计算流体力学领域得到应用[2][24]。但该算法求解强非线性动力学问题时产生非物理解（也称假解），严重影响求解精度。如图 1 所示，HDHB 法计算结果虽然是方程组(8)的数学解（使求解算法收敛），但已偏离了真实的物理响应。

Liu 指出[4]假解现象产生于对非线性项近似的过程中，存在高阶谐波对低次的混淆（污染）。以立方项非线性 $x^3$ 为例，一阶 HB 法中 $\hat{f}(\hat{x})$ 表达式
$$\hat{f}_{\text{HB}} = [r_0\ r_1\ r_2\ r_3\ r_4]^T$$

HDHB 法计算中，非线性项的频域分量(9)中不仅包含了原 HB 法中的所有项，还包括附加杂项，这些杂项的混入导致求解中出现非物理解。

$$\hat{\boldsymbol{f}}_{\mathrm{HDHB}} = \begin{bmatrix} r_0 + \frac{3}{4}\left(x_3^2 - x_4^2\right)x_1 - \frac{3}{2}x_2x_3x_4 \\ r_1 + \frac{3}{4}\left(x_1^2 - x_2^2\right)x_3 + \frac{3}{2}\left(x_3^2 - x_4^2\right)x_0 + \frac{1}{4}x_3^3 - \frac{3}{4}x_3x_4^2 - \frac{3}{2}x_1x_2x_4 \\ r_2 - \frac{3}{4}\left(x_1^2 - x_2^2\right)x_4 - \frac{1}{4}x_4^3 + \frac{3}{4}x_3^2x_4 - \frac{3}{2}x_1x_2x_3 - 3x_0x_3x_4 \\ r_3 + \frac{1}{4}x_1^3 + \frac{3}{4}\left(x_3^2 - x_4^2\right)x_1 - \frac{3}{4}x_1x_2^2 + \frac{3}{2}x_2x_3x_4 + 3x_0x_1x_3 - 3x_0x_2x_4 \\ r_4 + \frac{1}{4}x_2^3 + \frac{3}{4}\left(x_3^2 - x_4^2\right)x_2 - \frac{3}{4}x_1^2x_2 - \frac{3}{2}x_1x_3x_4 - 3x_0x_2x_3 - 3x_0x_1x_4 \end{bmatrix} \quad (9)$$

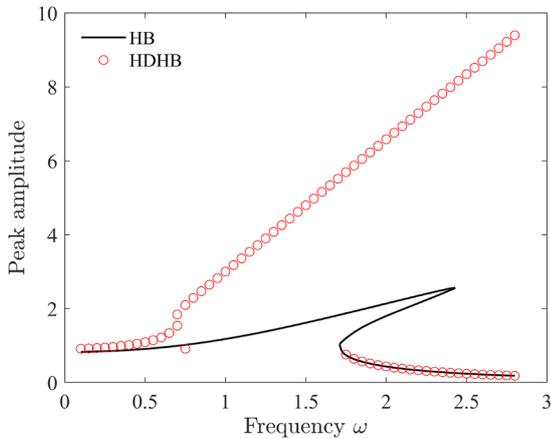

图 1 HB 法和 HDHB 法计算杜芬方程对比（$N=2$）

Fig.1 Comparison of the Duffing equation computed by the HB and HDHB methods（$N=2$）

### 1.3 重构谐波平衡法

Dai 证明[25]，HDHB 法与时域配点法等价，时域配点法通过使用时域配点上函数值对微分方程进行离散，从而联立代数方程组

$$\omega \boldsymbol{E}\boldsymbol{A}\hat{\boldsymbol{x}} = \tilde{\boldsymbol{f}}(\tilde{\boldsymbol{x}}) \quad (10)$$

其中配点矩阵

$$\boldsymbol{E} = \begin{bmatrix} 1 & \cos\omega t_0 & \sin\omega t_0 & \cdots & \cos N\omega t_0 & \sin N\omega t_0 \\ 1 & \cos\omega t_1 & \sin\omega t_1 & \cdots & \cos N\omega t_1 & \sin N\omega t_1 \\ \vdots & \vdots & \vdots & & \vdots & \vdots \\ 1 & \cos\omega t_M & \sin\omega t_M & \cdots & \cos N\omega t_M & \sin N\omega t_M \end{bmatrix}$$

混淆是造成计算出现非物理解的原因。至于高次谐波在时域配点计算中影响低次谐波的机理，则遵循"混淆规则"[25]：

**混淆规则**：假设 $\alpha \in [0,\pi]$ 被以 $h$ 为间隔均分为离散时间点，那么配点法中可区分的最大谐波次数 $n \in [-L,L]$，其中 $L = \pi/h$，$L$ 为"极限波次"。高次谐波 $n\left(|n| > L\right)$ 被误认为是相应低次谐波 $n_a$：

$$n_a = n - 2mL$$

其中 $n_a \in [-L,L]$，$m$ 是整数。

混淆规则指出，时域配点法可区分的最大谐波次数由配点数（离散时间点的间隔）决定。增加时域配点法中的配点数量，可区分的谐波次数越高，此时方程(10)个数多于待求解变量数，配点矩阵 $\boldsymbol{E}$ 为列满秩矩阵，存在 Moore-Penrose 广义逆矩阵 $\boldsymbol{E}^*$

$$\boldsymbol{E}^* = \frac{2}{M}\begin{bmatrix} 1/2 & 1/2 & \cdots & 1/2 \\ \cos\omega t_1 & \cos\omega t_2 & \cdots & \cos\omega t_M \\ \sin\omega t_1 & \sin\omega t_2 & \cdots & \sin\omega t_M \\ \vdots & \vdots & & \vdots \\ \sin N\omega t_1 & \sin N\omega t_2 & \cdots & \sin N\omega t_M \end{bmatrix}_{(2N+1)\times M}$$

使得

$$\boldsymbol{E}^*\boldsymbol{E} = \boldsymbol{I}_{2N+1}$$

其中 $\boldsymbol{I}_{2N+1}$ 为单位矩阵。在配点法方程同时作用矩阵 $\boldsymbol{E}^*$，对方程组降维，得到 RHB 法代数方程组

$$\omega \boldsymbol{A}\hat{\boldsymbol{x}} = \boldsymbol{E}^*\tilde{\boldsymbol{f}}(\tilde{\boldsymbol{x}}) \quad (11)$$

RHB 法在保证算法效率的同时，消除非物理解，从而实现对原 HB 法的最佳重构。**定理 1** 围绕如何选择合适的配点数，给出消除混淆确定条件[6]。

**定理 1**（条件等价性）：如果配点数 $M$，方法阶数 $N$，多项式非线性的幂次 $\phi$ 满足

$$M > (\phi+1)N \quad (12)$$

则 RHB 法与 HB 法等价。

为说明 RHB 法消除混淆的效果。分别使用 HDHB 法和 RHB 法计算杜芬方程，在分岔处（频率 $\omega = 2$）进行蒙特卡洛法模拟，选取 $10^4$ 组随机初

值（各频域分量在区间[-2,2]中选取），得到如表1所示的统计结果。统计结果表明，HDHB 法计算产生 58 组解，其中 55 组都是非物理解；而 RHB 法只计算得到 3 个具有物理含义的解[6]。

此外，HB-AFT 法与 HDHB 法计算流程一致，但 HB-AFT 法的原理是，选用配点数 $2\phi N+1$ 来消除混淆[11]，但过采样会占用计算机资源，在实际计算时会导致更大的 CPU 和 RAM 计算负担[2]。RHB 法基于时域配点法的统一框架，根据配点数的差异将所有 HB 类方法（HDHB 法、HB-AFT 法等）建立起联系。

表1 RHB 法与 HDHB 法在分岔处的解分布

Table 1 Statistical distribution of solution by the RHB and the HDHB method at bifurcation point

| Method | Upper branch(%) | Lower branch(%) | Unstable branch(%) | Non-physical(%) |
| --- | --- | --- | --- | --- |
| RHB | 57.13 | 19.71 | 23.16 | 0 |
| HDHB | 29.14 | 14.21 | 15.96 | 48.69 |

## 2 重构谐波平衡法改进策略

### 2.1 微分方程重铸技术

针对非多项式型非线性系统的 HB 法求解，Cochelin 等[12]提出将原系统改写为二次型系统：

$$\dot{z} = c + l(z) + q(z,z) \quad (13)$$

其中未知向量 $z$ 包含微分方程的原始变量 $x$ 及一些新的变量（引入的新变量用以改写方程）。$c$ 是关于未知量 $z$ 的常数向量；$l$ 是关于 $z$ 的线性向量值运算符；$q$ 则是二次向量值运算符。

因为对任意 $x^\phi$ 次非线性，RHB 法都能实现高效计算，改写后方程可以是更高次多项式，将式(13)进一步写成适配于 RHB 法的微分方程重铸形式：

$$\dot{z} = c + l(z) + n(z) \quad (14)$$

其中 $n$ 可以是关于变量 $z$ 的任意 $\phi$ 次多项式函数。重铸格式(14)涵盖多类型非多项式型非线性问题，下面将分类加以介绍。

（1）微分项耦合型

对范德波尔方程(15)分析，一阶导数 $\dot{x}$ 与平方项 $x^2$ 耦合。

$$\ddot{x} - \varepsilon(1-x^2)\dot{x} + x = F\cos\omega t \quad (15)$$

通过将方程重铸为标准形式(14)，得到

$$\begin{cases} \dot{x} = u \\ \dot{u} = (\varepsilon u - x) - \varepsilon x^2 u + F\cos\omega t \end{cases}$$

原方程转化为典型的非线性度 $\phi = 3$ 的动力学系统。得出结论：导数项的耦合不影响非线性度计入，在实际 HB 计算中，任意阶导数只计入一个非线性度。

（2）有理分式型

以 Rayleigh-Plesset 方程(16)为例

$$R\ddot{R} = -\frac{3}{2}\dot{R}^2 - A\frac{\dot{R}}{R} - B\frac{1}{R} + C\frac{1}{R^3} + D - E\cos(\omega t) \quad (16)$$

$A$，$B$，$C$，$D$ 均为常系数。将方程两端同除以 $R$，并运用倒数关系，令 $v = \dot{R}$，$x = 1/R$ 得到[6][12]

$$\begin{cases} \dot{R} = v \\ \dot{v} = -\frac{3}{2}v^2 x - Avx^2 - Bx^2 + Cx^4 + Dx - Ex\cos(\omega t) \\ 0 = Rx - 1 \end{cases}$$

处理有理分式型非线性项，通过额外增加方程 $Rx = 1$，将方程改写为非线性度 $\phi = 4$ 的多项式型非线性系统。

（3）根式型

相对论谐振子方程(17)中非线性项为根式

$$\ddot{x} + (1-\dot{x}^2)^{3/2} x = 0 \quad (17)$$

引入新变量 $u^2 = 1-\dot{x}^2$ 改写根号内关系式，利用根式的特性 $(1-\dot{x}^2)^{3/2} = (u^2)^{3/2} = u^3$，将方程改写为非线性度 $\phi = 4$ 的多项式型非线性系统。

$$\begin{cases} \dot{x} - v = 0 \\ \dot{v} + u^3 x = 0 \\ u^2 + v^2 - 1 = 0 \end{cases} \quad (18)$$

对任意形如 $\alpha^{q/p}$ 的根式，令 $\beta^p = \alpha$，使原根式型非线性转化得到多项式型：$\alpha^{q/p} = (\beta^p)^{q/p} = \beta^q$。

（4）初等超越函数

对初等超越函数的处理，需要导数信息来实现对原微分方程的重铸。以非线性单摆方程(19)为例

$$\ddot{\theta}(t) + \sin\theta(t) = 0 \quad (19)$$

由于三角函数是典型超越函数，不能通过简单的代数关系式改写转化为多项式型。首先额外引入两个变量 $s,c$ 来分别表示 $s(t) = \sin\theta(t)$，$c(t) = \cos\theta(t)$。利用三角函数的导数性质

$$\frac{\mathrm{d}}{\mathrm{d}t}\sin\theta(t)=\cos\theta(t)\cdot\dot{\theta}(t)$$
$$\frac{\mathrm{d}}{\mathrm{d}t}\cos\theta(t)=-\sin\theta(t)\cdot\dot{\theta}(t)$$

建立补充方程，从而实现对微分方程组的重铸[13]

$$\begin{cases}\dot{\theta}(t)=v(t)\\ \dot{v}(t)=-s(t)\\ \dot{s}(t)=c(t)v(t)\\ \dot{c}(t)=-s(t)v(t)\end{cases} \quad (20)$$

初等超越函数非线性的重铸，以导数关系作为方程组改写的依据（部分需要用到有理分式、根式型非线性重铸技巧）。附录表格中罗列了几类常见初等超越函数改写思路与重铸形式，可供参考。

## 2.2 多谐波平衡计算

"补频"（supplemental frequency）思想[17]是在原来单频假设的基础上，补充多个与响应相关的频率。假设考虑两个基频，将假设函数改写为

$$x_i(t)=\sum_m\sum_n x_{ic}(m,n)\cos((m\omega_1+n\omega_2)t)+x_{is}(m,n)\sin((m\omega_1+n\omega_2)t) \quad (21)$$

参数 $m$ 和 $n$ 满足不等式

$$|m|+|n|\leqslant p$$

其中 $p$ 代表 2 维傅里叶级数的截断[19][21]，类似于 RHB 法中阶数 $N$。

以 RHB 法为基础，提出 RMHB 法利用时域信息，将动力学问题(1)转化为非线性代数方程(11)求解。多基频计算中，配点矩阵 $\boldsymbol{E}$ 以及转换矩阵 $\boldsymbol{E}^*$ 写作 (23) 形式，记 $c^{m,n}=\cos(m\omega_1+n\omega_2)t$，$s^{m,n}=\sin(m\omega_1+n\omega_2)t$。在 RMHB 法计算中，存在选取合适的采样周期和时域配点数用以消除混淆的结论[23]。

**定理 2（多频计算中条件等价性）**：假设多项式非线性的幂次 $\phi$ 的非线性系统响应中包含两个基频，基频之比 $\omega_1/\omega_2$ 为有理数。则 RMHB 法消除混淆需满足采样时间 $T=2\pi/\mathrm{GCD}(\omega_1,\omega_2)$，配点数

$$M>(\phi+1)\frac{p\cdot\max(\omega_1,\omega_2)}{\mathrm{GCD}(\omega_1,\omega_2)} \quad (22)$$

其中 GCD 为两数最大公因数。

$$\begin{cases}\boldsymbol{E}=\begin{bmatrix}1 & 1 & \cdots & 1\\ c^{1,0}(t_1) & c^{1,0}(t_2) & \cdots & c^{1,0}(t_M)\\ s^{1,0}(t_1) & s^{1,0}(t_2) & \cdots & s^{1,0}(t_M)\\ \vdots & \vdots & \cdots & \vdots\\ c^{m,n}(t_1) & c^{m,n}(t_2) & \cdots & c^{m,n}(t_M)\\ s^{m,n}(t_1) & s^{m,n}(t_2) & \cdots & s^{m,n}(t_M)\\ \vdots & \vdots & \cdots & \vdots\\ c^{0,p}(t_1) & c^{0,p}(t_2) & \cdots & c^{0,p}(t_M)\\ s^{0,p}(t_1) & s^{0,p}(t_2) & \cdots & s^{0,p}(t_M)\end{bmatrix}^{\mathrm{T}}\\ \boldsymbol{E}^*=\frac{2}{M}\begin{bmatrix}\frac{1}{2} & \frac{1}{2} & \cdots & \frac{1}{2}\\ c^{1,0}(t_1) & c^{1,0}(t_2) & \cdots & c^{1,0}(t_M)\\ s^{1,0}(t_1) & s^{1,0}(t_2) & \cdots & s^{1,0}(t_M)\\ \vdots & \vdots & \cdots & \vdots\\ c^{m,n}(t_1) & c^{m,n}(t_2) & \cdots & c^{m,n}(t_M)\\ s^{m,n}(t_1) & s^{m,n}(t_2) & \cdots & s^{m,n}(t_M)\\ \vdots & \vdots & \cdots & \vdots\\ c^{0,p}(t_1) & c^{0,p}(t_2) & \cdots & c^{0,p}(t_M)\\ s^{0,p}(t_1) & s^{0,p}(t_2) & \cdots & s^{0,p}(t_M)\end{bmatrix}\end{cases} \quad (23)$$

## 3 数值仿真结果

### 3.1 相对论谐振子

作为物理学中经典问题，有必要对谐振子模型进行完整而严格的相对论处理[26]。考虑一个静质量为 $m$ 的质点在一维谐振力 $F=-k\bar{x}$ 的作用下做相对论运动。其中 $k$ 为弹性常数，$\bar{x}$ 为位移量。结合牛顿运动运动学方程以及动量定理，可以推导得到相对论振子方程[27-28]

$$\frac{\mathrm{d}^2\bar{x}}{\mathrm{d}\bar{t}^2}+\frac{k}{m}\left[1-\frac{1}{c^2}\left(\frac{\mathrm{d}\bar{x}}{\mathrm{d}\bar{t}}\right)^2\right]^{3/2}x=0$$

其中 $\bar{t}$ 是时间坐标（维度变量），对方程进行无量纲化得到方程(17)。初值条件 $x(0)=0$，$\dot{x}(0)=\beta$，其中 $-1<\beta<1$，且非线性振子的振动周期与对应周期解依赖于初始速度 $\beta$。Mickens 为了使相轨迹充满整个相平面引入非线性变换[29]

$$y=\frac{w}{\sqrt{1+w^2}}$$

再使用 HB 法计算。该处理方式在谐波振子低速运动时提供了高精度解[29-30]。当振子运动高速运动接近光速时（$\beta=0.85$），如图 2 所示，按照 Mickens 变换法求解将产生较大计算偏差。

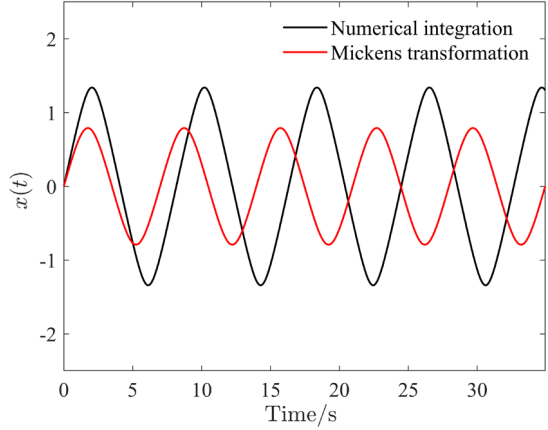

图 2 Mickens 变换解与真实物理解对比

Fig.2 Comparison of Mickens transformation result and real physical solution

RHB-重铸法按照根式型非线性重铸原则,将相对论谐振子改写为方程(18),再使用高阶 RHB 法求解。结合表 2 和图 3,RHB-重铸法可对高速运动的谐振子($\beta = 0.85$)直接进行高阶计算。五十五阶 RHB-重铸法能将计算误差控制在 $10^{-12}$ 量级(数值积分为参照),总计算耗时在 1 秒内。

表 2 各阶 RHB-重铸法求解相对论谐振子

Table 3 Solving relativistic harmonic oscillator by the RHB-recast method with different orders

| Order of method | Amplitude error | Computing time/s |
|---|---|---|
| 25 | $4.36\times10^{-7}$ | 0.73 |
| 35 | $1.95\times10^{-9}$ | 0.77 |
| 55 | $3.87\times10^{-12}$ | 1.07 |

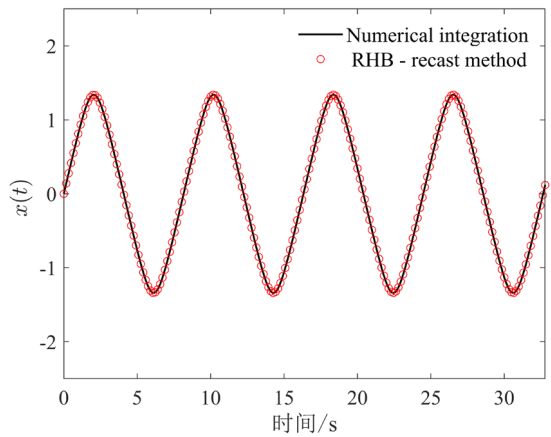

(a)时域响应

(a) Time domain response

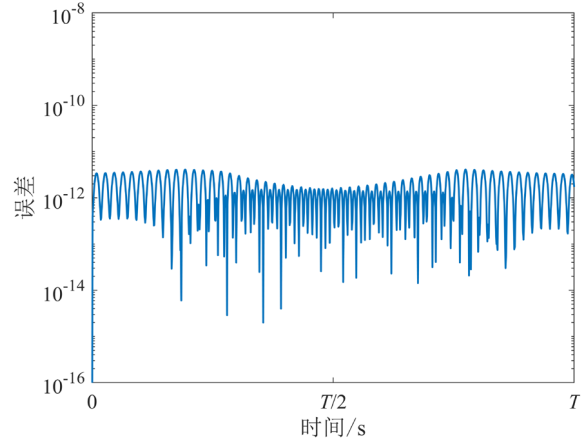

(b)误差曲线(一个周期内)

(b) Error curve(within one period)

图 3 五十五阶 RHB-重铸法求解相对论谐振子的计算结果

Fig.3 The computing result of the RHB-recast method ($N$=55) for solving relativistic harmonic oscillator

除通过方法阶数对求解精度的提高,还能通过两种方式实现对计算性能的改善:

(i)非线性代数方程组求解算法

非线性代数方程(11)求解算法能提高 HB 类方法计算性能,Zheng 等[31]结合 Tikhonov 正则化求解,黄建亮等[32]引入狗腿法思想结合回溯线搜索算法求解,Thomas 和 Dowell[33]使用 Broyden 法以提高 HB 法计算性能。本文根据算例来说明 L-M 法(Levenberg-Marquardt)较之传统迭代方法在求解 RHB 方程组时体现出优势。图 4 展示当 $\beta = 0.85$,频域量初值估计 $x_2$=1,$u_0$=0.7,使用两种不同的方程组求解算法:牛顿迭代法(Newton-Raphson)与 L-M 法得到的一个周期内误差曲线。不同于牛顿迭代法,L-M 法迭代格式为[34]

$$\bm{x}_{k+1} = \bm{x}_k - \left(\bm{J}_k^{\mathrm{T}}\bm{J}_k + \lambda_k \mathrm{diag}\left(\bm{J}_k^{\mathrm{T}}\bm{J}_k\right)\right)^{-1}\bm{J}_k^{\mathrm{T}}\bm{F}_k$$

其中 $\bm{x}_k$ 为当前迭代未知变量值,$\bm{F}_k$ 为函数值,$\bm{J}_k$ 为雅可比矩阵。L-M 法通过衡量每步迭代的误差是否发散,来决定撤回迭代并将参数 $\lambda_k$ 按十倍放缩。L-M 法在求解非线性方程组时显示出更高的计算精度,同比牛顿迭代法,精度提高了 $10^4$ 以上。

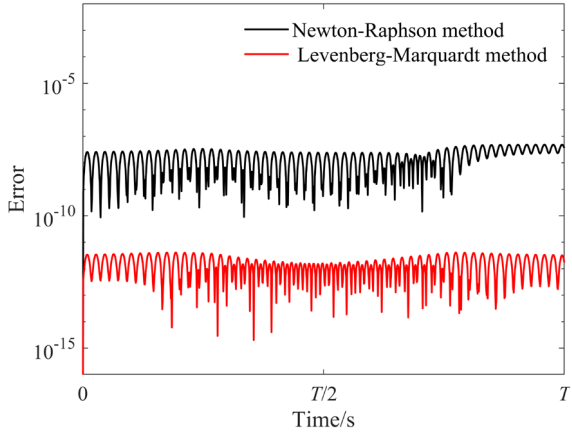

图 4 非线性方程组求解算法对计算精度的影响（相对论谐振子）

Fig.4 Influence of nonlinear equation algorithm on calculation accuracy（relativistic harmonic oscillator）

（ii）合理选择代数方程

相对论谐振子(17)为保守系统，振动频率 $\omega$ 是状态变量[35]。使用 RHB-重构法对 $n$ 自由度保守系统求解时，共需考虑 $n(2N+1)+1$ 个变量，需额外增加一个代数方程使 RHB 方程组适定。本文以初值约束来探讨方程对计算性能的影响。图 5 为采用五十五阶 RHB-重构法与 L-M 法求解器，使用不同代数方程（约束条件）得到的误差曲线。分别对应：初始位移约束 $x(0)=0$；初始速度约束 $\dot{x}(0)=\beta$；同时考虑初始位移与速度约束。同时考虑位移与速度约束时，RHB-重构法计算精度更高。

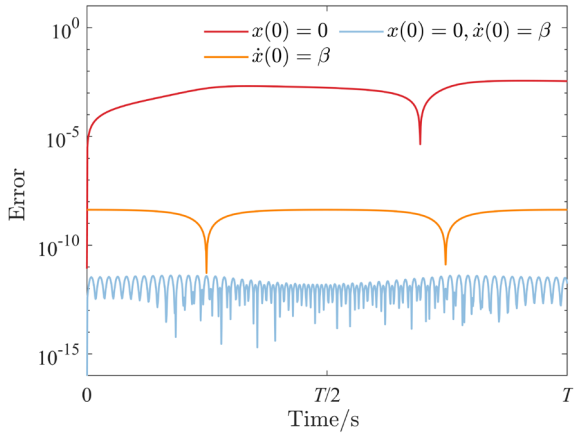

图 5 改变代数方程对计算精度影响

Fig.5 Influence of changed algebraic equation on calculation accuracy

相对论谐振子在高速运动时显示出复杂的动力学特性，需要高阶（>20 阶）方法来进行分析。使用重铸技术，将根式非线性转化为多项式型，再使用 RHB 法实现快速解算，克服高阶估计的限制。

如图 6 所示，RHB-重铸法适用于不同初速度条件下相对论振子的计算。

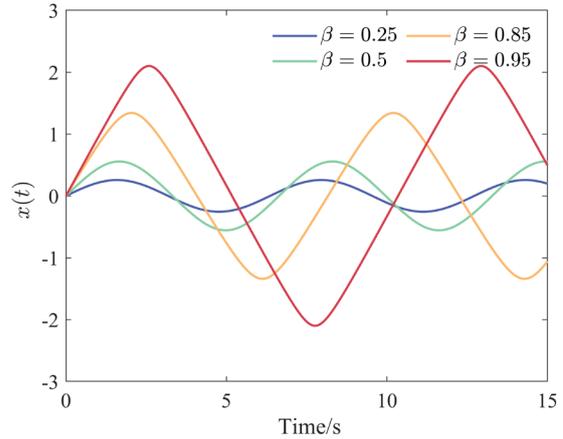

（a）位移-时间

（a）Displacement versus time

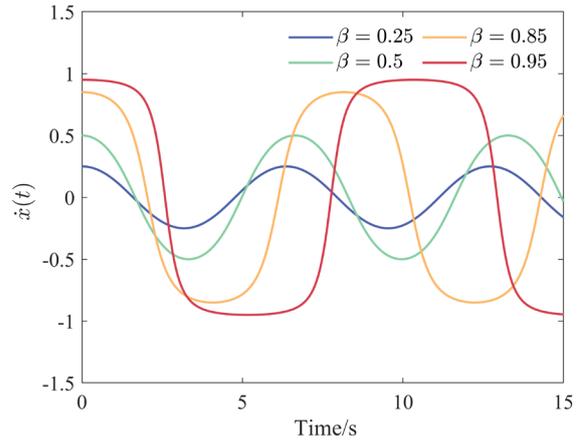

（b）速度-时间

（b）Velocity versus time

图 6 不同初速度条件下 RHB-重铸法求解相对论谐振子

Fig.6 Solving relativistic harmonic oscillator by the RHB-recast method under different initial velocities

### 3.2 非线性单摆

非线性单摆(19)是物理学中经典问题，且时域响应具有典型的周期性。但现有 HB 类方法不能做到高性能求解：HB-Taylor 法需要采用高的截断阶来保证计算精度；HB-AFT 法则需进行过采样以抑制混淆误差（如图 7 所示）。

对非线性单摆重铸，得到方程(20)。使用二十五阶 RHB-重铸法求解（大角度摆动，$\theta(0)=1.5$），初值估计 $x_1=1.423$，$s_1=1.065$，$c_0=1.028$，辅助 L-M 法求解器，计算结果如图 8 所示，与数值积分参考解比较，误差控制在 $10^{-12}$ 数量级以下。

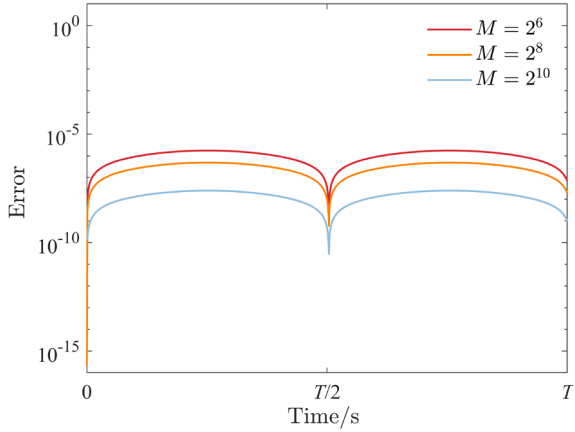

图 7 增加时域配点对抑制混淆误差的效果（使用 HB-AFT 法）

Fig.7 The effect of increasing time domain collocation on suppressing aliasing error（the HB-AFT method）

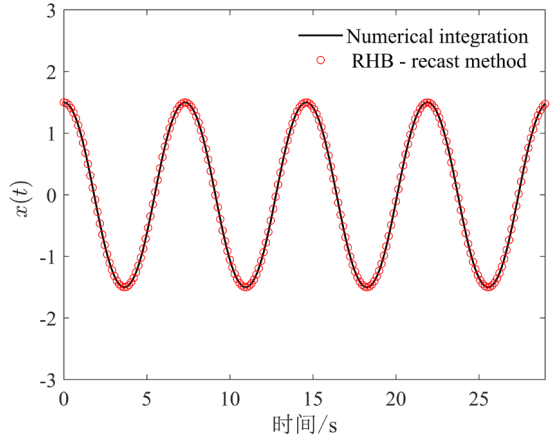

（a）时域响应

（a）Time domain response

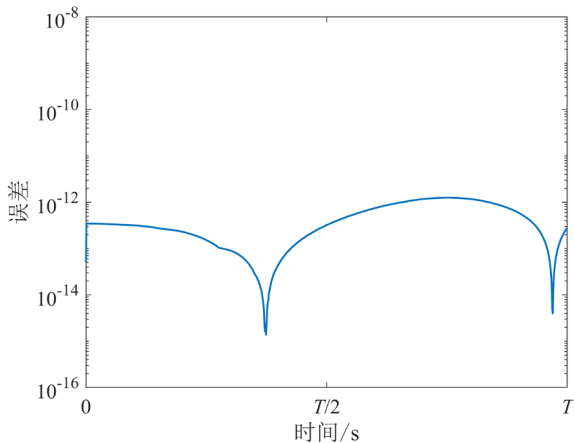

（b）误差曲线（一个周期内）

（b）Error curve（within one period）

图 8 二十五阶 RHB-重铸法求解非线性单摆

Fig.8 The computing result of the RHB-recast method（$N=25$）for solving nonlinear pendulum

分别采用牛顿迭代法、Tikhonov 正则法与 L-M 法得到图 9 计算结果。牛顿迭代法中雅可比阵奇异，求解失效；Tikhonov 法与 L-M 法计算所得的误差数量级相近，但是由于 Tikhonov 法对正则参数的优化使求解更耗时[31]。本例中，使用 L-M 法计算 RHB 法方程组仅需 0.72s，而 Tikhonov 法耗时达 1.17s，L-M 法较之 Tikhonov 法效率提高了 62.5%，达到了相近的求解精度。对比传统方法（牛顿迭代法）与优化方法（Tikhonov 法），L-M 法都更适合于 RHB 法方程组的求解。

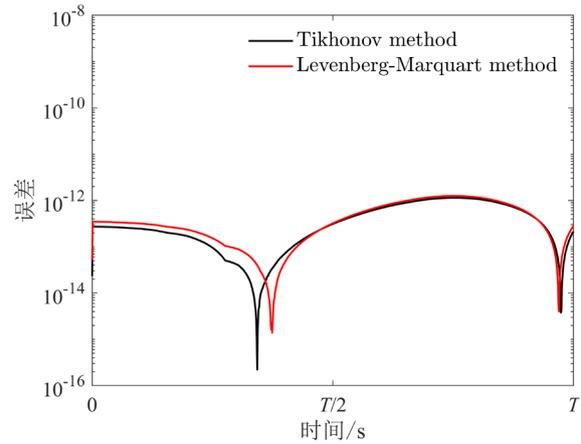

图 9 非线性方程组求解算法对计算精度的影响（非线性单摆）

Fig.9 Influence of nonlinear equation algorithm on calculation accuracy（nonlinear pendulum）

非线性单摆问题是保守系统，其周期解 $\theta(t)$ 与响应频率都依赖于初始振幅 $\theta(0)$。使用 RHB 法求解重铸方程组(20)时，本文给出三种代数方程组合方案，以助于对比分析：

**方案 1：** 对未知量 $\theta$, $v$, $s$ 和 $c$ 全谐波平衡（从常数项到 $N$ 次谐波项），计 $4(2N+1)$ 多个方程。增加初始振幅约束 $\theta(0) = \alpha$。

**方案 2：** 对未知量 $\theta$ 和 $v$ 全谐波平衡，计 $2(2N+1)$ 多个方程。$s$ 和 $c$ 从 1 次谐波开始，平衡系数到 $N$ 次，计 $2N$ 多个方程。增加初始振幅约束 $\theta(0) = \alpha$ 与两个对附加变量 $s$，$c$ 的初值约束：

$$\begin{cases} s(0) = \sin\theta(0) = \sin\alpha \\ c(0) = \cos\theta(0) = \cos\alpha \end{cases} \quad (24)$$

**方案 3：** 对未知量 $\theta$ 和 $v$ 全谐波展开，计 $2(2N+1)$ 多个方程。$s$ 和 $c$ 从 1 次谐波展开到 $N$ 次，计 $2N$ 多个方程。增加初始速度约束 $v(0) = 0$ 与两个对附加变量 $s$，$c$ 的初值约束(24)。

分别采用三种方案得到的计算结果如图 10 所

示。对比**方案 2** 和**方案 3**，**方案 1** 的计算误差大，主要原因是没有对附加变量 $s$ 和 $c$ 的初值进行约束。**方案 3** 同时利用了初始位移与速度信息，比**方案 2** 计算精度更高。

降低计算效率，计算用时达 1.80 秒。RHB-重铸法与 HB-AFT 法计算效率相当，但将计算精度提高了两个数量级以上。

RHB-重铸法适用于求解初等超越函数非线性问题。尤其是当单摆做大幅度振动时，传统的线性化手段无法很好地捕捉动力学响应，而 RHB-重铸法则可以提供高精度解析解（图 12）。

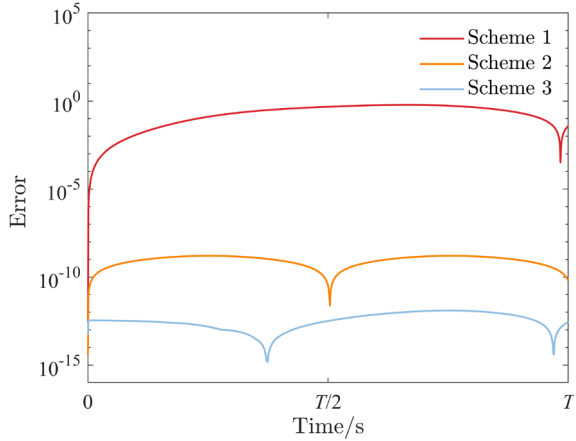

图 10 不同代数方程组合方案对计算精度影响
Fig.10 Influence of different combination schemes of nonlinear algebraic equation on calculation accuracy

对比求解非线性单摆的三种方法（HB-AFT 法、RHB-Taylor 法、RHB-重铸法）计算结果。考虑同阶截断 $N$=25，综合图 11 误差曲线和表 3 的计算结果，RHB-重铸法确定最优时域配点数 $M$=76，降低采样成本。RHB-Taylor 法虽可采用高阶截断（15 次多项式）保证计算精度，但超越函数的级数表示

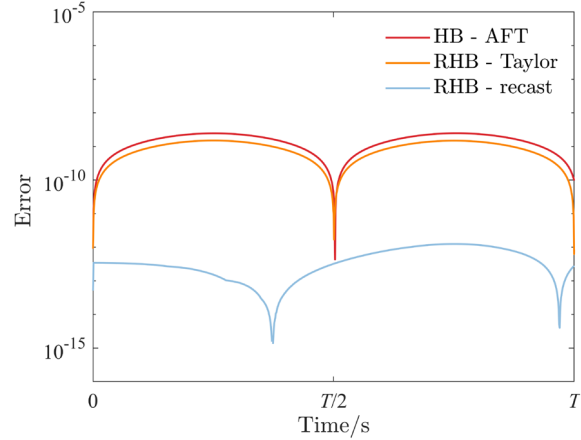

图 11 三种方法求解非线性单摆对应误差曲线对比
Fig.11 Comparison of corresponding error curves of three methods for solving nonlinear pendulum

表 3  三种方法计算非线性单摆结果对比
Table 3  Comparison of corresponding results of three methods for solving nonlinear pendulum

| Methods | $M$ | Average error | Computing time/s |
| --- | --- | --- | --- |
| RHB-recast | 76 | $1.9\times10^{-12}$ | 0.70 |
| RHB-Taylor | 401 | $8.9\times10^{-10}$ | 1.82 |
| HB-AFT | $2^{10}$ | $1.5\times10^{-9}$ | 0.69 |

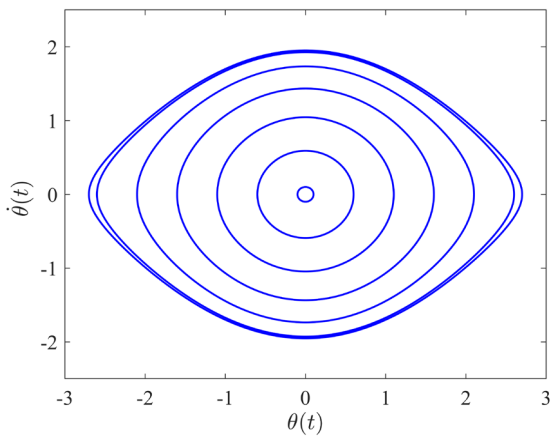

图 12 RHB-重铸法求解非线性单摆问题：相平面图

Fig.12 The RHB-recast method for solving nonlinear pendulum：phase plot

### 3.3 非线性耦合的非对称摆

用绳将一个质点小球悬挂在固定于弹性杆的末端，而弹性杆能在水平面与质点小球同步振动。此物理模型在许多其他二维摆系统中普遍存在，如球面摆、傅科摆等。特别地，当傅科摆由于不完美的悬挂或由于侧向运动引起非线性耦合而导致不对称，可能会导致额外的旋转，从而掩盖地面效应。非线性耦合的非对称单摆微分形式可以写作[36]

$$\begin{cases} \ddot{x}+(1-\kappa)x = -(1-\kappa)\left(x\dot{x}^2+x\dot{y}^2+xy\ddot{y}+x^2\ddot{x}\right) \\ \ddot{y}+y = -\left(y\dot{y}^2+y\dot{x}^2+xy\ddot{x}+y^2\ddot{y}\right) \end{cases} \quad (25)$$

此方程中二阶导数耦合入非线性项，导致使用传统方法进行直接求解变得棘手。以传统的有限差分类方法（欧拉法，龙格库塔法，MATLAB 内置 ode 算法等）而言，需额外计算代数方程组来辅助求解。而部分解析求解法的计算效果同样不佳，Jia 等[36]曾采用多时间尺度法推导方程(25)的解为

$$x = a\cos\sigma_x t + \frac{(1-\kappa)ab^2}{2b^2 - a^2 - 4\kappa}\cos(2\sigma_y - \sigma_x)t$$
$$y = b\cos\sigma_y t + \frac{(1-\kappa)a^2 b}{2a^2 - b^2 + 4\kappa}\cos(2\sigma_x - \sigma_y)t \quad (26)$$

令参数 $\kappa = 0.01$，摆的初始条件 $x(0)=0.1$，$y(0)=0.2$，$\dot{x}(0) = \dot{y}(0) = 0$。以修正 Chebyshev-Picard 迭代法（modified Chebyshev-Picard iteration, MCPI）计算结果为标准解[37-38]。仿真得到如图 13 所示，解(26)不仅推导复杂，且与真实的物理解间有较大误差。

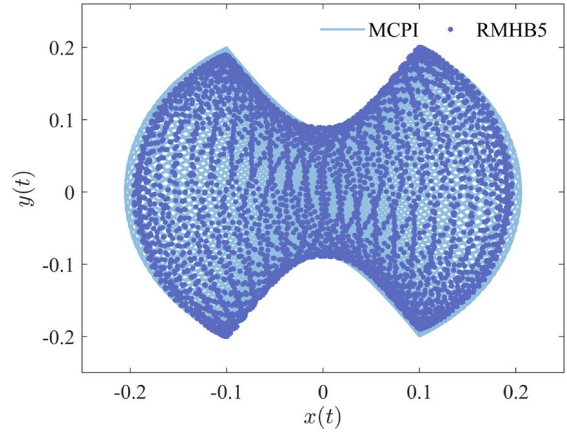

（a）运动轨迹
（a）Trajectory diagram

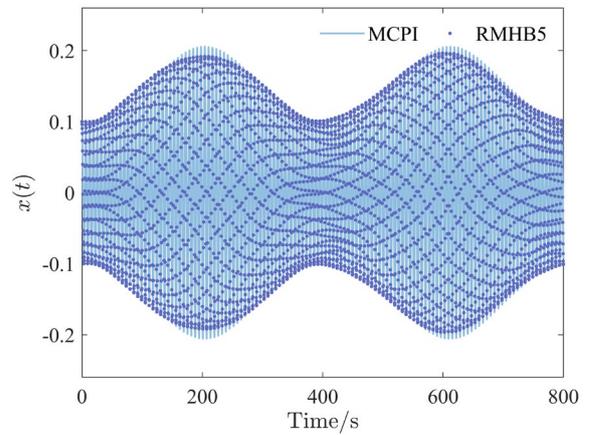

（b）x 方向振荡
（b）Oscillations in the x directions.

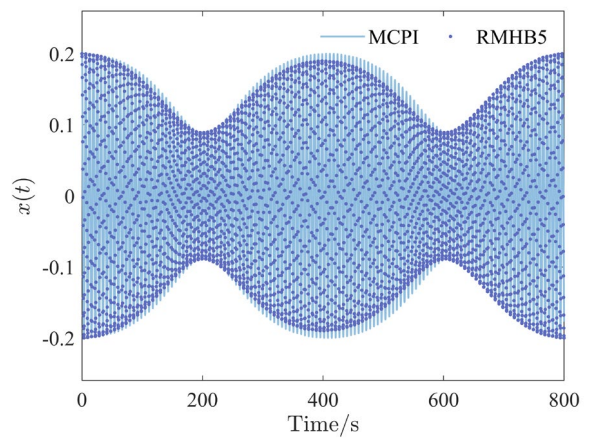

（c）y 方向振荡
（c）Oscillations in the y directions.

图 14 五阶 RMHB 法求解非线性耦合非对称摆
Fig.14 Solving nonlinear coupling asymmetric pendulum by the RMHB5

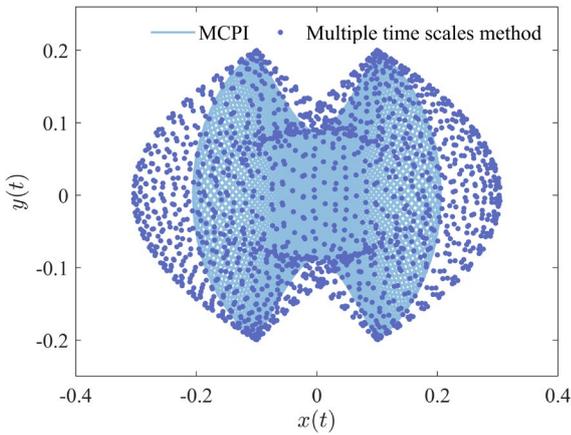

图 13 多时间尺度法与修正 Chebyshev-Picard 迭代法（参考解）求解轨迹对比
Fig.13 Comparison of the method of multiple scales and the MCPI method（reference solution）for solving trajectory diagram

对多时间尺度解(26)分析得知，动力学响应包含两个基频（所有频率成分都是基频的线性组合），通过快速傅里叶变换（FFT）得到：$\omega_1 = 0.9857$，$\omega_2 = 0.9935$。又从微分方程重铸的角度分析非对称摆微分方程(25)，因为高阶导数只算作一个非线性度，该方程是具有三次多项式非线性的动力学系统。使用五阶 RMHB 法进行计算得到图 14 所示仿真结果，计算用时 8.24 秒，$x$ 方向振动幅值误差 $9.21\times10^{-3}$，$y$ 方向振动幅值误差 $3.90\times10^{-7}$。

对二阶导耦合型非线性系统，有限差分方法不能直接求解，一些解析求解方法的计算精度低。考虑两个基频的 RMHB 法不仅保证计算效率，还能对拟周期运动进行准确的预测。

# 4 结 论

围绕复杂非线性动力学系统高效高精度周期解算需求，本文基于微分方程重铸、"补频"等思想，分别提出 RHB-重铸法与 RMHB 法，主要的工作与结论总结如下：

（1）提出 RHB-重铸法，将一般非线性问题改写为多项式型非线性方程，配合 RHB 法确定消除非物理解所需的最优配点阈值，实现对非多项式型非线性动力学系统周期响应进行高阶预测。数值实验说明，对比 HB-AFT 法的时域过采样、HB-Taylor 法对非线性函数进行泰勒级数截断两种处理方式，RHB-重铸法的计算精度高达 $10^{-12}$，且计算时间不超过 1 秒，同时保证了求解的高精度与高效率。

（2）提出结合"补频"思想的 RMHB 法，在原来单频假设的基础上，补充多个与物理响应相关的频率。借鉴 RHB 法中时域重构等价性，给出多频谐波平衡计算的最优配点数结论，并充分利用非线性量的时域采样代替复杂的高维傅里叶分析，最终实现对拟周期响应的准确捕捉。

（3）通过解算相对论谐振子、非线性单摆以及非线性耦合的非对称摆问题验证了本文算法的有效性。针对传统有限差分类方法求解困难的多自由度耦合系统，RHB 及本文提出的两种方法都不受方程形式的限制。此外，在相同的系统参数设置下，合理地选择代数方程求解器、代数约束方案，可有效提高非线性问题的求解精度。

本文提出的 RHB-重铸法与 RMHB 法对复杂非线性动力学系统周期性及拟周期响应的计算效率与精度相较于现行的方法与处理方式都具有优势。在后续的研究中，将探 RHB 法在高精周期响应求解方面的工程应用。

## 参考文献

附录

附表　常见初等超越函数的重铸

Recast form of the most common elementary transcendental functions

| Original function | Differential relationship | Companion variable | Recast equation |
| --- | --- | --- | --- |
| $u(t) = \exp[x(t)]$ | $\dot{u}(t) = \dot{x}(t) \cdot u(t)$ | | $\dot{u}(t) = \dot{x}(t) \cdot u(t)$ |
| $u(t) = \log_a[x(t)]$ | $\dot{u}(t) = \dot{x}(t)/(x(t) \cdot \ln a)$ | $v(t) = 1/x(t)$ | $\begin{cases} \dot{u}(t) = (v(t) \cdot \dot{x}(t))/\ln a \\ x(t)v(t) - 1 = 0 \end{cases}$ |
| $\begin{cases} u(t) = \sin[x(t)] \\ v(t) = \cos[x(t)] \end{cases}$ | $\begin{cases} \dot{u}(t) = \dot{x}(t) \cdot \cos[x(t)] \\ \dot{v}(t) = -\dot{x}(t) \cdot \sin[x(t)] \end{cases}$ | | $\begin{cases} \dot{u}(t) = \dot{x}(t) \cdot v(t) \\ \dot{v}(t) = -\dot{x}(t) \cdot u(t) \end{cases}$ |
| $u(t) = \tan[x(t)]$ | $\dot{u}(t) = \dot{x}(t) \cdot (1 + u^2(t))$ | | $\dot{u}(t) = \dot{x}(t) + \dot{x}(t) \cdot u^2(t)$ |
| $u(t) = \arcsin[x(t)]$ | $\dot{u}(t) = \dot{x}(t)/\sqrt{1 - x^2(t)}$ | $\begin{cases} v(t) = \sqrt{1 - x^2(t)} \\ w(t) = 1/v(t) \end{cases}$ | $\begin{cases} \dot{u}(t) = \dot{x}(t)w(t) \\ w(t)v(t) - 1 = 0 \\ v^2(t) + x^2(t) - 1 = 0 \end{cases}$ |
| $u(t) = \arccos[x(t)]$ | $\dot{u}(t) = -\dot{x}(t)/\sqrt{1 - x^2(t)}$ | $\begin{cases} v(t) = \sqrt{1 - x^2(t)} \\ w(t) = -1/v(t) \end{cases}$ | $\begin{cases} \dot{u}(t) = \dot{x}(t)w(t) \\ w(t)v(t) + 1 = 0 \\ v^2(t) + x^2(t) - 1 = 0 \end{cases}$ |
| $u(t) = \arctan[x(t)]$ | $\dot{u}(t) = \dot{x}(t)/(1 + x^2(t))$ | $\begin{cases} v(t) = 1 + x^2(t) \\ w(t) = 1/v(t) \end{cases}$ | $\begin{cases} \dot{u}(t) = \dot{x}(t)w(t) \\ w(t)v(t) - 1 = 0 \\ v(t) - x^2(t) - 1 = 0 \end{cases}$ |